\documentclass[12pt,onecolumn,oneside,leqno,a4paper]{article}
\usepackage{german}
\usepackage{exscale}
\usepackage{amsmath,amssymb,verbatim}
\usepackage{amscd}
\usepackage{graphicx}
\usepackage{rotating}
\usepackage{amsxtra,amsthm}
\usepackage{color}
\usepackage{makeidx}
\usepackage[pdftex]{hyperref}
\hypersetup{bookmarksopen=true,bookmarksopenlevel=1,bookmarksnumbered=true}

\newcommand{\TB}[1]{\textbf{#1}}
\newcommand{\TI}[1]{\textit{#1}}

\def\Bew{\textit{Proof. }}

\newcommand\QED{\hspace*{\fill} $ \square $}
\newcommand\euro{{\sffamily C%
  \makebox[0pt][l]{\kern-.70em\mbox{--}}%
  \makebox[0pt][l]{\kern-.68em\raisebox{.25ex}{--}}}}

\newenvironment{DES1}[1]
	{\begin{list}{}
	{\settowidth{\labelwidth}{#1}
	\setlength{\labelsep}{1em} 
	\setlength{\leftmargin}{0em}
	\setlength{\itemindent}{\labelwidth}
	\addtolength{\itemindent}{\labelsep}
	\setlength{\parsep}{\parskip}
	\setlength{\itemsep}{0.2\parsep plus0.02\parsep minus0.02\parsep}
	\setlength{\topsep}{0.2\parsep plus0.02\parsep minus0.02\parsep}
	\setlength{\partopsep}{0.6\parsep plus0.02\parsep minus0.05\parsep}
	}}
	{\end{list}}

\def\R{\mathbb{R}}

\def\Rplus{\mathbb{R}^+}

\def\B{\mathbb{B}}

\def\SS{\mathbb{S}}
\renewcommand{\H}{\mathbb{H}}
\newcommand{\Rz}{\R^{2}}
\newcommand{\Rd}{\R^{3}}

\def\const{\text{\textup{const.}}}

\newcommand{\pfeil}[1]{\overrightarrow{#1}}

\newcommand\Iff{\;\iff\;}

\newcommand\Arctan{\arctan\,}

\newcommand\Arccos{\arccos\,}

\newcommand\Sin{\sin\,}
\newcommand\Cos{\cos\,}

\newcommand\vhi{\varphi}
\newcommand\vro{\varrho}
\newcommand\Sinh{\sinh\,}
\newcommand\Cosh{\cosh\,}
\newcommand\Tanh{\tanh\,}

\renewcommand\Dot[2]{\left\langle #1,#2 \right\rangle}

\def\ds{\displaystyle}

\newcommand{\Part}[2]{\frac{\partial #1}{\partial #2}}

\renewcommand{\Hat}[1]{\widehat{\,#1\,}}

\newcommand{\cross}{\times}

\newcommand{\Ast}{{\ds \ast}}

\newcommand{\Ac}{\mathcal{A}}

\newcommand{\downto}{\downarrow }
\newcommand{\upto}{\uparrow }
\newcommand{\norm}[1]{\left\|#1\right\|}
\newcommand{\BM}{\begin{pmatrix}}
\newcommand{\EM}{\end{pmatrix}}
\newcommand{\BV}{\begin{vmatrix}}
\newcommand{\EV}{\end{vmatrix}}

\newtheoremstyle{tld}
{}
{}
{\itshape}
{}
{\bfseries}
{.}
{2em}
{}

\newtheoremstyle{bb}
{}
{}
{}
{}
{\bfseries}
{.}
{2em}
{}

\newtheoremstyle{bbs}
{}
{}
{}
{}
{\bfseries}
{.}
{\newline}
{}

\swapnumbers

\theoremstyle{tld}
\newtheorem{thm}{\negthickspace. Theorem}[section] 
\newtheorem{thmz}[thm]{\negthickspace. Zusatz}
\newtheorem{lem}[thm]{\negthickspace. Lemma}
\newtheorem{cor}[thm]{\negthickspace. Corollary}
\newtheorem{defi}[thm]{\negthickspace. Definition}
\newtheorem{defthm}[thm]{\negthickspace. Definition und Satz}
\newtheorem{thmdef}[thm]{\negthickspace. Satz und Definition}

\theoremstyle{bb}
\newtheorem{ex}[thm]{\negthickspace. Example}
\newtheorem{rem}[thm]{\negthickspace. Remark}
\newtheorem{aufg}{\negthickspace. Aufgabe}[section] 

\theoremstyle{bbs}
\newtheorem{exs}[thm]{\negthickspace. Examples}
\newtheorem{rems}[thm]{\negthickspace. Remarks}

\newcommand{\BT}{\begin{thm}}
\newcommand{\ET}{\end{thm}}

\newcommand{\BTZ}{\begin{thmz}}
\newcommand{\ETZ}{\end{thmz}}

\newcommand{\BL}{\begin{lem}}
\newcommand{\EL}{\end{lem}}

\newcommand{\BC}{\begin{cor}}
\newcommand{\EC}{\end{cor}}

\newcommand{\BD}{\begin{defi}}
\newcommand{\ED}{\end{defi}}

\newcommand{\BDT}{\begin{defthm}}
\newcommand{\EDT}{\end{defthm}}

\newcommand{\BTD}{\begin{thmdef}}
\newcommand{\ETD}{\end{thmdef}}

\newcommand{\BX}{\begin{ex}}
\newcommand{\EX}{\end{ex}}

\newcommand{\BXS}{\begin{exs}}
\newcommand{\EXS}{\end{exs}}

\newcommand{\BR}{\begin{rem}}
\newcommand{\ER}{\end{rem}}

\newcommand{\BRS}{\begin{rems}}
\newcommand{\ERS}{\end{rems}}

\newcommand{\BE}{\begin{equation}}
\newcommand{\EE}{\end{equation}}

\newcommand{\BA}{\begin{aufg}}
\newcommand{\EA}{\end{aufg}}

\numberwithin{equation}{section}

\newcommand{\bild}[1]{\includegraphics{BILDER/#1}}

\newcommand{\abs}[1]{\left|#1\right|}
\newcommand{\del}{\partial}

\begin{document} 

\thispagestyle{empty}

\begin{center}
\TB{\LARGE Polygons in hyperbolic geometry 2:}

\smallskip
\TB{\LARGE Maximality of area}

\TB{\large Rolf Walter}
\end{center}

\setcounter{section}{0}


\section{\hspace{-1em}.
Introduction}
\label{intro}


The topic of this part is the maximum question for the area
of polygons in the hyperbolic plane with \TI{fixed}
sidelengths. As a main result it will be shown: Among all
polygons in the hyperbolic plane with fixed positive
sidelengths there exist polygons of maximal area. Each such
maximal polygon is either oriented-convex and cocyclic or
else collinear. In the first case the maximal area is
positive, in the second case it is zero. 
A more detailed version will be given in Theorem 
\ref{thm_polymax_main}. As a corollary
one obtains that among the non-collinear polygons the only
copies for which a rigidity can be hoped for are the
oriented-convex cocyclic ones.
The cocyclicity means that the vertices are situated on a
distance circle, a distance line or a horocycle, the three
types of cycles in the hyperbolic plane. The phenomenon of
different circle types stands in salient contrast to the
Euclidean case and pays for various difficulties in the 
hyperbolic discussion.

The corresponding result in the Euclidean plane
has been discussed several times in the literature:
In \TI{Yaglom/Boltjanski} [1951] a proof
of the cocyclicity for maximal polygons in $ \Rz $ is given
within the class of simply connected polygons, using the
general isoperimetric inequality. It is stated there that a
proof without this tool would be extremely difficult. 
Other treatments for the Euclidean
case are in \TI{Blaschke} [1956], \TI{Kryzhanovskij} [1959], and
\TI{Knebelman} [1941]. Of course, there are some ideas from the
Euclidean situation which are also worthwhile in the
non-Euclidean case, but new phenomena and difficulties occur.
In particular, this is true for the notion of area itself.
Also familiar constructions from Euclidean geometry are no
more available. For instance there is no circumferential
angle theorem in the hyperbolic circle theory and no
similarities exist for figures in the hyperbolic space.

In classical expositions of hyperbolic geometry, there
prevails the relation of area to the angle sum. However,
this relation is directly applicable only for polygons which
are bounding because it rests on the Gauss/Bonnet integral
theorem. The polygons to be considered here are more
general: no a-priori assumptions on their form (convexity,
simple closedness, etc.) have to be made. Therefore, the
angle sum is hard to handle for our purpose. Since the
sidelengths are strongly involved it is more advisable to
view the whole problem in the context
of distance geometry in the sense of \TI{Menger} [1928] and 
\TI{Blumenthal} [1970]. So an effort is necessary to express the area by
different means, in particular by distances instead of
angles. This is achieved by a more analytical definition of the 
area functional in combination with explicit expressions to 
be found in a paper of \TI{Bilinski} [1969].

For basics on hyperbolic geometry we refer to part 1. 
In particular
the circle model of Cayley/Klein will be considered
throughout. It has the advantage that geodesics are
Euclidean lines such that hyperbolic convexity properties are 
very near to their Euclidean relatives. This is not true for
the angles but they are almost not entering anyway.


\section{\hspace{-1em}.
The polygon area}
\label{polarea}
\markright{\ref{polarea}. The polygon area}

The Cayley/Klein model of the hyperbolic plane consists of
the open unit ball $ \B $ in $ \Rz $ where the hyperbolic lines 
are just the chords of the horizon $ \SS := \del \B $ (which
itself doesn't belong the hyperbolic plane). In this model,
the riemannian metric on $ \B $ is  given by
$$
ds^{2} = 
\frac{1}{(1-\xi^{2}-\eta^{2})^{2}}
\Bigl(
(1-\eta^{2})d\xi^{2}+2\xi\eta d\xi
d\eta+(1-\xi^{2})d\eta^{2}
\Bigr),
$$
where $ (\xi,\eta) $ are cartesian coordinates in $ \B $. 
The area element of this metric, viewed as a $ 2 $-form, sounds
$$
\mu = \frac{d\xi \wedge d\eta}{\sqrt{1-\xi^2-\eta^2}^3}.
$$
It is the exterior derivative of a certain $ 1 $-form,
namely
\begin{subequations}
\label{eq_polymax_leib}
\BE
\mu = d\omega, \qquad 
\omega := 
\frac{-\eta\;d\xi+\xi\;d\eta}{1-\xi^{2}-\eta^{2}+\sqrt{1-\xi^{2}-\eta^{2}}}.
\EE
This is the decisive access to the area here. The form $ \omega $ 
plays the same role in the hyperbolic plane as the form 
$ -\eta\;d\xi+\xi\;d\eta $ does in the Euclidean plane. It 
admits the calculation of the area $ F(A) $ 
of any compact subset
$ A \subset \B $ with `good' oriented boundary $ \del A = K $ by a
curve integral:
\BE
\label{eq_polymax_leib1}
F(A) = \int_{K} \omega,
\EE
\end{subequations}
in complete analogy to the Leibniz formula in the Euclidean 
case. A `good' boundary is e.g. a closed continuous $ C^{1} $-chain
without selfintersections, in particular a polygon chain without 
selfintersections (see \TI{Cartan} [1967], Sects. 4.2--4.4).

However, the integral in \eqref{eq_polymax_leib} is much more
general. It yields, for any oriented chain $ K $ composed of compact 
$ C^{1} $-arcs with real weights, a real value which depends
linearly under the addition and scalar multiplication of
such chains. This value is well defined insofar as it is
independent of the representation of the chain. In this
generality, the value can also be zero or negative,
depending on the orientation of the chain. So we are dealing
with \TI{signed} areas.

By Bilinski [1969], Eqn. (6.2), the signed area of a
$ 3 $-gon (i.e. a triangle) in $ \B $ with vertices 
$ A,B,C $ is explicitly given by
\BE
\label{eq_polyex_bil1}
F(ABC) = 2\,\Arctan
\frac{[a,b,c]}{\Dot{a}{b}+\Dot{b}{c}+\Dot{c}{a}+1},
\EE
where $ a,b,c \in \Rd $ are the normalized point vectors of the
vertices. For the pseudo-Euclidean scalar product involved
here and its rules see part 1, Sect. 2.

Now, combining the chain integral from \eqref{eq_polymax_leib1}
with the expression \eqref{eq_polyex_bil1}, one 
obtains the following explicit representation for the signed 
area of any $ n $-gon $ P := Z_{1}\ldots Z_{n} $:
\BE
\label{eq_polyex_pinh}
F(P) =
2\sum_{k=1}^{n} \Arctan
\frac{[z,z_{k},z_{k+1}]}%
{
\Dot{z}{z_{k}}+
\Dot{z}{z_{k+1}}+
\Dot{z_{k}}{z_{k+1}}+1
}.
\EE
Here, $ Z $ can be any point in $ \B $, and all point
vectors occurring must be normalized. The properties of the
chain integral \eqref{eq_polymax_leib1} ensure that this
expression is indeed independent of the choice of the
`origin' $ Z $. One may call Eqn. \eqref{eq_polyex_pinh} 
the \TI{parachute formula} since all connecting triangles
with $ Z $ are summed together with the right account of
signs. The corresponding Euclidean formula is much simpler
and is sometimes named after Gau"s.

This independence immediately yields the following facts:

\medskip
\BL
\label{lem_polymax_nn}
~\\[0.4ex]
Each oriented-convex polygon has positive area. 
Each collinear polygon has area zero.
\EL

\Bew
The special choice $ Z := Z_{1} $ in \eqref{eq_polyex_pinh} 
leads to
$$
F(P) =
2\sum_{k=2}^{n-1} \Arctan
\frac{[z_{1},z_{k},z_{k+1}]}%
{
\Dot{z_{1}}{z_{k}}+
\Dot{z_{1}}{z_{k+1}}+
\Dot{z_{k}}{z_{k+1}}+1
}.
$$
In the first case, every numerator in the sum is
positive, so is $ F(P) $. In the second case every numerator
is $ 0 $, so is $ F(P) $.
\QED

\medskip
Also, the existence of polygons with maximal area can be
deduced from \eqref{eq_polyex_pinh}, using the general
maximum principle:

\medskip
\BL
\label{lem_polymax_genmax}
For each $ n $-gon there exists a $ n $-gon 
of maximal area with same sidelengths.
\EL

\Bew
A $ n $-gon  $ Z_{1}\ldots Z_{n} $ can be represented 
by a point in the cartesian
product $ \B^{n} := \B \cross \cdots \cross \B $ ($ n $
factors), say equipped with the maximum metric $ d_{n} $.
Without loss of generality it is possible to fix the point 
$ Z_{1} $ for all polygons to be considered. As is obvious
from Eqn. \eqref{eq_polyex_pinh}, the area 
$ F(Z_{1},\ldots,Z_{n}) $ depends continuously on 
$ (Z_{2},\ldots,Z_{n}) $. The definition set consists of all
points $ (Z_{2},\ldots,Z_{n}) \in\B^{n-1} $ with 
$ d(Z_{k},Z_{k+1}) = L_{k} = \const $, $ k = 1,\ldots,n $.
It is compact, namely bounded and closed: The boundedness
follows from the estimate
$$
\begin{aligned}
d_{n-1}((Z_{2},\ldots,Z_{n}),(Z_{1},\ldots,Z_{1})) 
&= \max
\{
d(Z_{2},Z_{1}),
d(Z_{3},Z_{1}),
\ldots,
d(Z_{n},Z_{1})
\}\\[1ex]
&\leq
(n-1)\max \{L_{1},\ldots,L_{n}\},
\end{aligned}
$$
and the closedness is deduced from the continuity of the
functions 
$ (Z_{2},\ldots,Z_{n}) \mapsto d(Z_{k},Z_{k+1}) $, 
$ k = 1,\ldots,n $.
Thus $ F(Z_{1},\ldots,Z_{n}) $  has on this definition set
(which is not vacuous) a finite maximum.
\QED

\medskip
The main question is of course: How do the maximal 
$ n $-gons look like? As in the Euclidean case the final
answer will be: In general they are cocyclic. But on the way
to this goal one needs an analytical characterization for
the cocyclicity, at least for low $ n $. Non-collinear
triangles always have a circum-circle. So the next
interesting case are the quadrangles. For three and four
points in the hyperbolic plane there are several identities
and characterizations which will prepare the answer 
(see Sect. \ref{tripquad}).

\newpage
\BRS
\label{rem_polymax_huil}
~\\[-6ex]
\begin{DES1}{(ii)}
\item[(i)]
In Eqn. \eqref{eq_polyex_bil1}, the absolute value of the
determinant $ [a,b,c] $ is expressible by Gram's identity
in terms of the pairwise scalar products of the point vectors 
$ a,b,c $, so finally in terms of the sidelengths. With some
fancy conversions for hyperbolic functions this is converted 
to a formula of classical \TI{L'Huilier type}:
\begin{subequations}
\BE
\label{eq_polymax_huil1}
\abs{F(ABC)} =
4\,\Arctan\sqrt{
\Tanh\frac{S}{4}
\Tanh\frac{D_{1}}{4}
\Tanh\frac{D_{2}}{4}
\Tanh\frac{D_{3}}{4},
}
\EE
where
\BE
\label{eq_polymax_huil2}
\begin{aligned}
S &:= \phantom{-} L_{1}+L_{2}+L_{3} \\[1ex]
D_{1} &:= -L_{1}+L_{2}+L_{3} \\[1ex]
D_{2} &:= \phantom{-} L_{1}-L_{2}+L_{3} \\[1ex]
D_{3} &:= \phantom{-} L_{1}+L_{2}-L_{3},
\end{aligned}
\EE
\end{subequations}
thus expressing the absolute value of the area of a triangle 
solely by its sidelengths $ L_{1},L_{2},L_{3} $.

The quantities $ D_{1},D_{2},D_{3} $ are just the
differences from the triangle inequalities, so 
$ D_{1} \geq 0 $, $ D_{2} \geq 0 $, $ D_{3} \geq 0 $.
The appropriate definition set of the right hand side of 
\eqref{eq_polymax_huil1} is
$$ 
\Lambda := 
\{(L_{1},L_{2},L_{3}) \mid 
L_{1} > 0, \; L_{2} > 0, \; L_{3} > 0, \;
D_{1} \geq 0, \; D_{2} \geq 0, \; D_{3} \geq 0\}.
$$
Denote by $ H: \Lambda \to \R $ the 
\TI{L'Huilier-function}, given by the right hand side of 
\eqref{eq_polymax_huil1}. 
For $ D_{1} > 0 $, $ D_{2} > 0 $, $ D_{3} > 0 $, $ H $
depends real-holomorphically on
$ (L_{1},L_{2},L_{3}) $. At the boundary of $ \Lambda $ 
($ D_{1} = 0 $ or $ D_{2} = 0 $ or $ D_{3} = 0 $) $ H $ 
is still continuous but no more differentiable. But just
this singularity will be helpful later on for the growth of 
the area. The growth of $ H $ is controlled by the following
limit relations:
\begin{align}
\label{eq_polymax_huil4}
\lim_{X \upto L_{2}+L_{3}}\;
\Part{H}{X}(X,L_{2},L_{3}) &= -\infty
\quad\text{if}\quad
L_{2} > 0, \; L_{3} > 0\\[1ex]
\label{eq_polymax_huil3}
\lim_{X \downto L_{1}-L_{2}}\;
\Part{H}{X}L_{1},L_{2},X) &=\quad \infty 
\quad\text{if}\quad
L_{1} > L_{2} > 0.
\end{align}
They immediately follow from the corresponding partial derivatives 
in the interior of $ \Lambda $.
\item[(ii)]
Besides the L'Huilier expression, there are some other
formulas for the absolute value of the triangle area which
will be needed; see Bilinski [1969], Eqns. (11.3) and (11.5):
\begin{align}
\label{eq_polymax_bil_113}
\abs{F(ABC)} &= 2\Arctan\,
\frac{\Sinh L_{1} \Sinh L_{2} \,\Sin \gamma}%
{\Cosh L_{1}+\Cosh L_{2}+\Cosh L_{3}+1}, 
\\[1ex]
\label{eq_polymax_bil_115}
\abs{F(ABC)} &= 2\,\Arccos
\frac{\Cosh L_{1}+\Cosh L_{2}+\Cosh L_{3}+1}%
{\ds
4\,\Cosh\frac{L_{1}}{2}\Cosh\frac{L_{2}}{2}\Cosh\frac{L_{3}}{2}},
\end{align}
where $ \gamma \in \left]0,\pi\right[ $ is the angle opposite 
to the side of length $ L_{3} $. 
Instead of the `over-determined' formula 
\eqref{eq_polymax_bil_113} one may use a variant which
arises by substituting $ \Cosh L_{3} $ according to the 
cosine law:
\BE
\label{eq_polyex_bil4}
\abs{F(ABC)} = 2\,\Arctan
\frac{\Sinh L_{1}\Sinh L_{2}\Sin \gamma}%
{(1+\Cosh L_{1})(1+\Cosh L_{2})-\Sinh L_{1}\Sinh
L_{2}\Cos\gamma}.
\EE
Eqn. \eqref{eq_polyex_bil4} is the sole instance where
angles enter the game. It will be needed for
`parallelogram-like' quadrangles which occur as exceptional
cases in the area maximizing problem.
\end{DES1}
\ERS


\section{\hspace{-1em}.
Identities for triples und quadruples}
\label{tripquad}
\markright{\ref{tripquad}. Identities for triples und quadruples}

Here, certain properties of triples and quadruples of
points will be expressed solely by distances. In particular 
this applies to the property of cocyclicity. The following
abbreviations will be used:
\BE
\label{eq_tripquad_F0}
K_{PQ} := \Cosh d(P,Q), \qquad
S_{PQ}:=\Sinh \frac{d(P,Q)}{2},  \qquad P,Q \in \B.
\EE
In most cases, the hyperbolic distances $ d(P,Q) $ enter in
form of the $ \sinh $-quantities $ S_{PQ} $, but in a few cases the
expressions become easier with the $ \cosh $-quantities
$ K_{PQ} $.

\medskip
\TB{\TI{Point triples}}

The vertices of a triangle are collinear if their images 
on the quadric shell $ \H \subset \Rd $ lie in a vector plane. If the
vertices are not collinear then they are always cocyclic
because their images lie in a plane of $ \Rd $ which
doesn't contain $ 0 $. The image points in $ \Rd $ are also not
collinear since $ \H $ cuts any straight line of $ \Rd $ at 
most twice. So the plane and hence the circum-circle is 
uniquely determined. Certainly, this circle needs not to
have a center, i.e. it can also be a distance line or a
horocycle. Which case occurs can be read off solely from the
sidelengths.

\medskip
\BL
\label{lem_tripquad_idrei}
Let $ A,B,C $ be non-collinear points in $ \B $ and 
$ a,b,c $ their image points on $ \H $. Then the plane
in $ \Rd $, affinely spanned by $ a,b,c $ has the equation
\BE
\label{eq_tripquad_idrei1}
\Dot{u}{x} = p \qquad \text{with} \qquad
u := a\cross b+b\cross c+c\cross a, \quad
p := [a,b,c].
\EE
With the corresponding distance quantities
\eqref{eq_tripquad_F0}, the following representations hold
true:
\BE
\label{eq_tripquad_idrei2}
\begin{aligned}
\Dot{u}{u} &= 8
(S_{AB}^2 S_{BC}^2+S_{BC}^2 S_{CA}^2+S_{CA}^2 S_{AB}^2)-
4(S_{AB}^4+S_{BC}^4+S_{CA}^4) \\[1ex]
&= 4
(S_{AB}+S_{BC}+S_{CA})
(S_{AB}+S_{BC}-S_{CA})
(S_{BC}+S_{CA}-S_{AB})
(S_{CA}+S_{AB}-S_{BC})\\[1ex]
[a,b,c]^{2} &= 16S_{AB}^2 S_{BC}^2 S_{CA}^2+\Dot{u}{u}.
\end{aligned}
\EE
\EL

\Bew
The vector $ u $ from \eqref{eq_tripquad_idrei1} doesn't
vanish because
$$
a\cross b+b\cross c+c\cross a = (b-a) \cross (c-a).
$$
The plane affinely spanned by $ a,b,c $ has indeed 
the equation  $ \Dot{u}{x} = p $ with the values of 
$ u $ and $ p $ as in \eqref{eq_tripquad_idrei1} because 
each of $ a,b,c $ satisfies it, e.g.
$$
\Dot{a\cross b+b\cross c+c\cross a}{a} = 
\Dot{b\cross c}{a} =
[a,b,c].
$$
This proves the first part.

\smallskip
From $ u $ in \eqref{eq_tripquad_idrei1} follows for 
$ \Dot{u}{u} $ by the Lagrange identity, since 
$ \norm{a} = \norm{b} = \norm{c} = 1 $:
$$
\begin{aligned}
\Dot{u}{u} 
&=
1-\Dot{a}{b}^{2}+1-\Dot{b}{c}^{2}+1-\Dot{c}{a}^{2} \\[1ex]
&\phantom{=}+
2(\Dot{a}{b}\Dot{b}{c}-\Dot{c}{a})+
2(\Dot{b}{c}\Dot{c}{a}-\Dot{a}{b})+
2(\Dot{c}{a}\Dot{a}{b}-\Dot{b}{c}).
\end{aligned}
$$
Substituting the scalar products according to the formula
\BE
\label{eq_tripquad_skS}
\Dot{a}{b} = \Cosh d(A,B) = 
1+2\sinh^{2} \frac{d(A,B)}{2} = 1+2S_{AB}^{2}
\EE
and also its cyclic extensions then yields 
the first representation in \eqref{eq_tripquad_idrei2}
and by factorizing the second one.

Finally, $ [a,b,c]^{2} $ is by Gram's identity expressible
in terms of scalar products and then in the same manner by
the corresponding $ \sinh $-quantities. 
The third line follows from this by comparison with the
second line.
\QED

\medskip
In part 1, Theorem 4.3 it has been stated that the
sidelengths already determine on which circle type the
vertices of an oriented-convex cocyclic polygon are situated.
Now, for triangles, this can be made completely explicit:

\medskip
\BC
\label{cor_tripquad_fdrei}
Define for non-collinear points $ A,B,C $ in $ \B $ the invariant
$$
\Delta := 
(S_{AB}+S_{BC}-S_{CA})
(S_{BC}+S_{CA}-S_{AB})
(S_{CA}+S_{AB}-S_{BC}).
$$
Then the quantity $ \Delta $ which solely depends on the
pairwise distances of the points determines the type of the
circum-circle of $ A,B,C $, namely:
$$
\begin{aligned}
\Delta > 0 &\Iff \text{the circum-circle is a distance
circle} \\
\Delta < 0 &\Iff \text{the circum-circle is a distance line} \\
\Delta = 0 &\Iff \text{the circum-circle is a horocycle}.
\end{aligned}
$$
\EC

\Bew
Comparison of the middle equation in
\eqref{eq_tripquad_idrei2} with the classification of the
circle types in part 1, Sect. 2.
\QED

\medskip
\TB{\TI{Point quadruples}}

Let $ a,b,c,e $ \TI{points} in the pseudo-Euclidean space $ \Rd $.
The decision whether these points are \TI{coplanar}, i.e.
are contained in a plane, depends on the \TI{quadruple quantity}
\BE
\label{eq_tripquad_ivier1}
[a,b,c,e] := [a,b,c]-[b,c,e]+[c,e,a]-[e,a,b].
\EE
Namely, the vectors $ b-a,c-a,e-a $ are linearly dependent
iff the determinant $ [b-a,c-a,a-e] $ vanishes. Expanding
this determinant by the multilinear and alternating rules just 
results in the quadruple quantity. Thus:
\BE
\label{eq_tripquad_ivier4}
\text{$ a,b,c,e \in \Rd $ coplanar} \Iff [a,b,c,e] = 0.
\EE
On the other hand, the \TI{vectors} $ a,b,c,e $ are always
linearly dependent. This can be formulated by means of the hyperbolic  
analogue of the so called Cayley/Menger determinant (see 
Blumenthal [1970], Ch. IV, \S\;40 and Ch. XII, \S\;106). Since,
here, the relation to the quadruple quantity is needed, an
equivalent expression will be deduced directly as follows: A
`universal' dependency relation formulated without any
scalar product is $ V(a,b,c,e) = 0 $, where
\BE
\label{eq_tripquad_ivier2}
V(a,b,c,e) := [a,b,c]e-[b,c,e]a+[c,e,a]b-[e,a,b]c.
\EE

\BL[$ \mathbf{17} $-identity]
\label{lem_tripquad_R17}
For any points $ A,B,C,E $ in $ \B $ one always has:
$$
\begin{aligned}
0 &=K_{AB}^2 K_{CE}^2+K_{AC}^2 K_{BE}^2+K_{AE}^2 K_{BC}^2 \\[1ex]
  &\quad-2 K_{AB} K_{AC} K_{BE} K_{CE}-2 K_{AB} K_{AE} K_{BC} K_{CE}
   -2 K_{AC} K_{AE} K_{BC} K_{BE} \\[1ex]
  &\quad+2 K_{AB} K_{AC} K_{BC}+2 K_{AB} K_{AE} K_{BE}
   +2 K_{AC} K_{AE} K_{CE}+2 K_{BC} K_{BE} K_{CE} \\[1ex]
  &\quad-K_{AB}^2-K_{AC}^2-K_{AE}^2-K_{BC}^2-K_{BE}^2-K_{CE}^2 
  +1.
\end{aligned}
$$
\EL

\Bew
This is just the relation 
$ \Dot{V(a,b,c,e)}{V(a,b,c,e)} = 0 $, followed by expansion 
with means of Gram's identity, and expressing the occurring 
scalar products according to the first part of 
\eqref{eq_tripquad_skS}.
\QED

\medskip
The announced relation then sounds:

\medskip
\BL
\label{lem_tripquad_ivier}
Let $ A,B,C,E $ be points in $ \B $ and $ a,b,c,e $ their
images in $ \H $. Then between the quantities 
\eqref{eq_tripquad_ivier1} and \eqref{eq_tripquad_ivier2}
the following identity holds true:
\BE
\label{eq_tripquad_ivier3}
\begin{aligned}
&\Dot{V(a,b,c,e)}{V(a,b,c,e)} = \\[1ex]
&4(S_{AB}S_{CE}+S_{AC}S_{BE}+S_{AE}S_{BC})[a,b,c,e]^{2}+\\[1ex]
&64
(S_{AC}S_{BE}-S_{AE}S_{BC}-S_{AB}S_{CE})
(S_{AB}S_{CE}-S_{AC}S_{BE}-S_{AE}S_{BC})\cdot\\[1ex]
&(S_{AE}S_{BC}-S_{AC}S_{BE}-S_{AB}S_{CE}).
\end{aligned}
\EE
\EL

\Bew
If the computation of Lemma \ref{lem_tripquad_R17} is
continued by replacing the $ \cosh $-quantities by the 
$ \sinh $-quantities according to the second part of
\eqref{eq_tripquad_skS}, one arrives at
\BE
\label{eq_tripquad_gencond}
\begin{aligned}
&\frac{1}{64}\Dot{V(a,b,c,e)}{V(a,b,c,e)} = \\[1ex]
&S_{AB}^{4}S_{CE}^{4}-2S_{AB}^{2}S_{AC}^{2}S_{BE}^{2}S_{CE}^{2}
-2S_{AB}^{2}S_{AE}^{2}S_{BC}^{2}S_{CE}^{2}+S_{AC}^{4}S_{BE}^{4}
\\[1ex]
&-2S_{AC}^{2}S_{AE}^{2}S_{BC}^{2}S_{BE}^{2}+S_{AE}^{4}S_{BC}^{4}+S_{AB}^{4}S_{CE}^{2}
+S_{AB}^{2}S_{AC}^{2}S_{BC}^{2}
-S_{AB}^{2}S_{AC}^{2}S_{BE}^{2}
\\[1ex]
&-S_{AB}^{2}S_{AC}^{2}S_{CE}^{2}-S_{AB}^{2}S_{AE}^{2}S_{BC}^{2}+S_{AB}^{2}S_{AE}^{2}S_{BE}^{2}-S_{AB}^{2}S_{AE}^{2}S_{CE}^{2}
-S_{AB}^{2}S_{BE}^{2}S_{CE}^{2}
\\[1ex]
&+S_{AB}^{2}S_{CE}^{4}+S_{AC}^{4}S_{BE}^{2}-S_{AC}^{2}S_{AE}^{2}S_{BC}^{2}-S_{AB}^{2}S_{BC}^{2}S_{CE}^{2}-S_{AC}^{2}S_{AE}^{2}S_{BE}^{2}+S_{AC}^{2}S_{AE}^{2}S_{CE}^{2}
\\[1ex]
&-S_{AC}^{2}S_{BC}^{2}S_{BE}^{2}+S_{AC}^{2}S_{BE}^{4}+S_{AE}^{4}S_{BC}^{2}+S_{AE}^{2}S_{BC}^{4}-S_{AE}^{2}S_{BC}^{2}S_{BE}^{2}
-S_{AE}^{2}S_{BC}^{2}S_{CE}^{2}\\[1ex]
&+S_{BC}^{2}S_{BE}^{2}S_{CE}^{2}-S_{AC}^{2}S_{BE}^{2}S_{CE}^{2}.
\end{aligned}
\EE
The calculation of the square $ [a,b,c,e]^{2} $
follows the same pattern as before by expanding 
with Gram's identity and replacing the
scalar products according to \eqref{eq_tripquad_skS}. The
result is:
\BE
\label{eq_tripquad_planecond}
\begin{aligned}
&\frac{1}{16}[a,b,c,e]^{2} = \\[1ex]
&S_{AB}^{4}S_{CE}^{2}+S_{AB}^{2}S_{AC}^{2}S_{BC}^{2}
-S_{AB}^{2}S_{AC}^{2}S_{BE}^{2}-S_{AB}^{2}S_{AC}^{2}S_{CE}^{2}
-S_{AB}^{2}S_{AE}^{2}S_{BC}^{2}\\[1ex]
&+S_{AB}^{2}S_{AE}^{2}S_{BE}^{2}
-S_{AB}^{2}S_{AE}^{2}S_{CE}^{2}-S_{AB}^{2}S_{BC}^{2}S_{CE}^{2}
-S_{AB}^{2}S_{BE}^{2}S_{CE}^{2}+S_{AB}^{2}S_{z,w}^{4}\\[1ex]
&+S_{AC}^{4}S_{BE}^{2}-S_{AC}^{2}S_{AE}^{2}S_{BC}^{2}
-S_{AC}^{2}S_{AE}^{2}S_{BE}^{2}+S_{x,z}^{2}S_{AE}^{2}S_{CE}^{2}
-S_{AC}^{2}S_{BC}^{2}S_{BE}^{2}\\[1ex]
&+S_{AC}^{2}S_{BE}^{4}
-S_{AC}^{2}S_{BE}^{2}S_{CE}^{2}+S_{AE}^{4}S_{BC}^{2}
+S_{AE}^{2}S_{BC}^{4}-S_{AE}^{2}S_{BC}^{2}S_{BE}^{2}\\[1ex]
&-S_{AE}^{2}S_{BC}^{2}S_{CE}^{2}+S_{BC}^{2}S_{BE}^{2}S_{CE}^{2}.
\end{aligned}
\EE
Then, by direct calculation, the assertion 
\eqref{eq_tripquad_ivier3} is deduced from 
\eqref{eq_tripquad_gencond} and \eqref{eq_tripquad_planecond}.
\QED

\medskip
This implies a first characterization of the cocyclicity:

\medskip
\BC[Ptolemy relations]
\label{cor_tripquad_pt}
For any points $ A,B,C,E $ in $ \B $ define
\BE
\label{eq_tripquad_ptqs}
\begin{aligned}
\Delta_{1}&:=S_{AC}S_{BE}-S_{AE}S_{BC}-S_{AB}S_{CE} \\[1ex]
\Delta_{2}&:=S_{AB}S_{CE}-S_{AC}S_{BE}-S_{AE}S_{BC} \\[1ex]
\Delta_{3}&:=S_{AE}S_{BC}-S_{AC}S_{BE}-S_{AB}S_{CE}.
\end{aligned}
\EE
Then there hold the inequalities
\BE
\label{eq_tripquad_pt1}
\Delta_{1} \leq 0, \qquad 
\Delta_{2} \leq 0, \qquad
\Delta_{3} \leq 0.
\EE
If $ A,B,C,E $ are not collinear then $ A,B,C,E $ are 
cocyclic if and only if, in \eqref{eq_tripquad_pt1},
the equals sign occurs at least once, i.e.
\BE
\label{eq_tripquad_pt4}
\Delta_{1} = 0 \quad \text{or}\quad 
\Delta_{2} = 0 \quad \text{or}\quad
\Delta_{3} = 0.
\EE
\EC

\Bew

\TI{For \eqref{eq_tripquad_pt1}:}
Eqn. \eqref{eq_tripquad_ivier3} says indeed
\BE
\label{eq_tripquad_prpt}
-(S_{AB}S_{CE}+S_{AC}S_{BE}+S_{AE}S_{BC})[a,b,c,e]^{2} =
16\Delta_{1}\Delta_{2}\Delta_{3},
\EE
hence always 
\BE
\label{eq_tripquad_pt2}
\Delta_{1}\Delta_{2}\Delta_{3} \leq 0.
\EE
Moreover 
\BE
\label{eq_tripquad_pt3}
\begin{aligned}
\Delta_{1}+\Delta_{2} &= -2S_{AE}S_{BC} \leq 0 \\[1ex]
\Delta_{2}+\Delta_{3} &= -2S_{AC}S_{BE} \leq 0 \\[1ex]
\Delta_{1}+\Delta_{3} &= -2S_{AB}S_{CE} \leq 0.
\end{aligned}
\EE
\TI{Case $ \Delta_{1}\Delta_{2}\Delta_{3} < 0 $:}
If one of these factors were positive, say 
$ \Delta_{1} > 0 $ then, from \eqref{eq_tripquad_pt2}, 
the other two must be of different
sign, so another factor has to be positive, say 
$ \Delta_{2} > 0 $, Then, from \eqref{eq_tripquad_pt3}, a
contradiction can be read off. Hence the assertion in this
case, and indeed $ \Delta_{1} < 0 $, $ \Delta_{2} < 0 $, 
$ \Delta_{3} < 0 $.

\TI{Case $ \Delta_{1}\Delta_{2}\Delta_{3} = 0 $:}
At least one of these factors vanishes, say 
$ \Delta_{1} = 0 $. Then from \eqref{eq_tripquad_pt3}:
$ \Delta_{2} \leq 0 $, $ \Delta_{3} \leq 0 $.

\smallskip
\TI{For the remaining assertion:}

If the points $ A,B,C,E $ are cocyclic then the points 
$ a,b,c,e $ are coplanar in $ \Rd $, thus from 
\eqref{eq_tripquad_ivier4} and \eqref{eq_tripquad_prpt}: 
$ \Delta_{1}\Delta_{2}\Delta_{3} = 0 $.

For the converse, assume $ \Delta_{1}\Delta_{2}\Delta_{3} = 0 $
and consider the two cases for the affine hull 
$ \Ac $ of $ a,b,c,e $ in $ \Rd $ : 
(a) $ \dim \Ac = 1 $; (b) $ \dim \Ac \geq 2 $.

In case (a), there exist two different points among
$ a,b,c,e $, say $ a \neq b $ in $ \H $, and then
$ c,d \in a \vee b $. Since any straight line in $ \Rd $ 
cuts $ \H $ at most twice: $ c,e \in \{a,b\} $. So 
$ A,B,C,E $ are collinear: Case (a) cannot happen.

In case (b), there exist three points among $ a,b,c,e $ in
general position, say $ a,b,c $, a fortiori pairwise
different. Then
$ S_{AB}S_{CE}+S_{AC}S_{BE}+S_{AE}S_{BC} > 0 $ because 
not both of $ S_{CE}, S_{BE} $ can vanish. 
Thus, by \eqref{eq_tripquad_prpt}: $ [a,b,c,e] = 0 $.
The points $ a,b,c,e $ affinely span a plane in 
$ \Rd $ which doesn't contain $ 0 $ since $ A,B,C,E $ are
not collinear. So $ A,B,C,E $ are cocyclic.
\QED

\medskip
The Ptolemy equations \eqref{eq_tripquad_pt4}
have the disadvantage that each of them 
contains all \TI{six} pairwise distances while a
quadrangle is generally determined by \TI{five} distances.
The following results work against this disadvantage.

\medskip
\BC
\label{cor_tripquad_cinq}
If the non-collinear points $ A,B,C,E \in \B $ are
cocyclic then at least one of the following equations is
valid:
\begin{align}
\label{eq_tripquad_cinq1}
(S_{AB}S_{BC}+S_{AE}S_{CE})S_{AC}^2 &= 
(S_{AB}S_{CE}+S_{AE}S_{BC})(S_{AB}S_{AE}+S_{BC}S_{CE})
\\[1ex]
\label{eq_tripquad_cinq2}
-(S_{AB}S_{BC}-S_{AE}S_{CE})S_{AC}^2 &= 
(S_{AB}S_{CE}-S_{AE}S_{BC})(S_{AB}S_{AE}-S_{BC}S_{CE}).
\end{align}
Also, at least one of the two equations which arise from 
\eqref{eq_tripquad_cinq1}, \eqref{eq_tripquad_cinq2} by
permuting the points $ A,B,C,D $ holds true, in particular
at least one of the equations
\begin{align}
\label{eq_tripquad_Cinq1}
(S_{AB}S_{AE}+S_{BC}S_{CE})S_{BE}^2 &= 
(S_{AB}S_{CE}+S_{BC}S_{AE})(S_{AB}S_{BC}+S_{AE}S_{CE})
\\[1ex]
\label{eq_tripquad_Cinq2}
-(S_{AB}S_{AE}-S_{BC}S_{CE})S_{BE}^2 &= 
(S_{AB}S_{CE}-S_{BC}S_{AE})(S_{AB}S_{BC}-S_{AE}S_{CE}).
\end{align}
is valid.
\EC

\Bew
Clearly, it suffices to proof that  
\eqref{eq_tripquad_cinq1} or \eqref{eq_tripquad_cinq2} holds.

These equations 
arise by eliminating $ S_{BE} $ from $ [a,b,c,e] = 0 $ and
$ \Dot{V(a,b,c,e)}{V(a,b,c,e)} = 0 $, with the left hand sides 
expressed by 
\eqref{eq_tripquad_planecond} and \eqref{eq_tripquad_gencond}.
Both equations  are biquadratic w.r.t. variable $ S_{BE} $, 
i.e. they only contain $ S_{BE}^{4} $ and $ S_{BE}^{2} $.
The leading coefficients are $ S_{AC}^2 $ resp. 
$ (1+S_{AC}^2)S_{AC}^2 $. In case $ S_{AC} \neq 0 $, the
elimination can be done via the resultant of two quadratic
polynomials (cf. e.g. \TI{van der Waerden} [2003], \S\,30). 
The resultant comes out as
$$
\begin{aligned}
R &:= S_{AC}^4\cdot
((S_{AB}S_{BC}+S_{AE}S_{CE})S_{AC}^2
-(S_{AB}S_{CE}+S_{AE}S_{BC})(S_{AB}S_{AE}+S_{BC}S_{CE}))^2
\\[1ex]
&\phantom{:= S_{AC}^4\;}\cdot ((S_{AB}S_{BC}-S_{AE}S_{CE})S_{AC}^2
+(S_{AB}S_{CE}-S_{AE}S_{BC})(S_{AB}S_{AE}-S_{BC}S_{CE}))^2.
\end{aligned}
$$
This implies the assertion if $ S_{AC} \neq 0 $. 

In case $ S_{AC} = 0 $, i.e. $ C = A $, the points 
$ A, B, E $ are always situated on a circle (if not
collinear), and also Eqn. \eqref{eq_tripquad_cinq2} is
always satisfied.
\QED

\medskip
Assuming convexity, one can say more:

\medskip
\BL[$ \text{Perron [1964]} $]
~\\[0.3ex]
\label{lem_tripquad_perron}
For any oriented-convex cocyclic $ 4 $-gon $ ABCE $ in 
$ \B $ there hold Eqns. \eqref{eq_tripquad_cinq1}
\underline{and} \eqref{eq_tripquad_Cinq1}.
\EL

\Bew
Using the means of part 1, in particular Lemma 3.1, 
this can be done be straightforward calculation.
A cocyclic $ 4 $-gon $ ABCE $ is oriented-convex iff, in a
suitable representation of the circum-circle, the group
parameters of the vertices are in monotonic order. E.g. for
a distance circle of hyperbolic radius $ R $ in standard
position, the quantities $ S_{AB} $, etc. sound by
Eqn. (2.32) of part 1:
$$
S_{AB} = \vro \cdot \Sin \frac{\vhi_{B}-\vhi_{A}}{2},\ldots,
\qquad \vro := \Sinh R,
$$
where $ \vhi_{A},\ldots,\vhi_{E} $ are the parameter values 
of the points $ A,\ldots,E $ in the representation (2.22) 
of part 1. By $ \vhi_{A} < \vhi_{B} < \vhi_{C} < \vhi_{E} $ 
and $ 0 < \vhi_{E}-\vhi_{A} < 2\pi $, all these sine-values 
are positive and, by due trigonometric conversions, both
sides of \eqref{eq_tripquad_cinq1} resp. \eqref{eq_tripquad_Cinq1}
turn out to be equal, namely
$$
\vro^2\cdot \sin^{2} \frac{\vhi_{C}-\vhi_{A}}{2}
\quad\text{resp.}\quad
\vro^2\cdot \sin^{2} \frac{\vhi_{E}-\vhi_{B}}{2}.
$$
Similar calculations confirm Eqns.
\eqref{eq_tripquad_cinq1}, \eqref{eq_tripquad_Cinq1} for
distance lines and horocycles as circum-circles.
\QED

\medskip
Even more important is the converse since the equations of
Corollary \ref{cor_tripquad_cinq} contain one variable less 
then the Ptolemy equations of Corollary \ref{cor_tripquad_pt}:

\medskip
\BT
\label{thm_tripquad_umkperron}
Let $ ABCE $ be a $ 4 $-gon in $ \B $ with $ A \neq C $ such that 
$ B $ and $ E $ lie on different sides of the diagonal line
$ A \vee C $. Then, the relation
\BE
\label{eq_tripquad_umk1}
S_{AC}^2 = 
\frac{(S_{AB}S_{CE}+S_{AE}S_{BC})(S_{AB}S_{AE}+S_{BC}S_{CE})}%
{S_{AB}S_{BC}+S_{AE}S_{CE}}, 
\EE
implies that the $ 4 $-gon $ ABCE $ is cocyclic.
\ET

\Bew
By Theorem 5.3 of part 1, there exists an
oriented-convex cocyclic $ 4 $-gon $ A'B'C'E' $ with same
sidelengths as $ ABCE $. Its diagonal length $ d(A',C') $
is, by Lemma \ref{lem_tripquad_perron}, 
computed from the sidelengths 
$ d(A',B') $ , $ d(B',C') $, $ d(C',E') $, $ d(A',E') $
by the same formula as, by assumption \eqref{eq_tripquad_umk1}, 
$ d(A,C) $ is calculated from the sidelengths 
$ d(A,B) $ , $ d(B,C) $, $ d(C,E) $, $ d(A,E) $. 
This implies $ d(A',C') = d(A,C) $. So the triangles 
$ ABC $ und $ A'B'C' $ are congruent, and the same is true
for the triangles $ CEA $ und $ C'E'A' $.
Moreover, the points $ B' $, $ C' $ lie on different sides
of the line $ A' \vee C' $. For the vertex $ E $, there are 
left two possible positions which only differ by reflection 
on $ A \vee C $. Exactly one of this positions has the
property that $ E $, $ B $ are on different sides of the 
line $ A \vee C $. For this position the $ 4 $-gon $ ABCE $ 
is congruent to the $ 4 $-gon $ A'B'C'E' $. So, the 
$ 4 $-gon $ ABCE $ must be cocyclic.
\QED


\section{\hspace{-1em}.
Polygons with maximal area}
\label{polymax}
\markright{\ref{polymax}. Polygons with maximal area}

We already know from Lemma
\ref{lem_polymax_genmax} that the maximum
problem of the polygon area for fixed sidelengths is
solvable. Here, the maximal copies will be determined as
expressed in detail in the main result
\ref{thm_polymax_main}. When we speak of maximal polygons or
of the enlargement of polygons this always refers to the
area functional for fixed sidelengths.

\medskip
In general, the vertices of a $ n $-gon will be denoted by 
$ Z_{1},\ldots,Z_{n} $. In order to keep compliance with
the special annotations of Sect. \ref{tripquad}, we follow the
identifications $ A = Z_{1} $, $ B = Z_{2} $, $ C = Z_{3} $,
$ E = Z_{4} $ without further mention.

The following fact is helpful in order to exclude eventual 
degenerate cases:

\medskip
\BL
\label{lem_polymax_pwd}
A $ n $-gon for which two non-adjacent vertices coincide
can always be enlarged.
\EL

\Bew
Without loss of generality, assume $ Z_{1} = Z_{m} $ for an 
index $  m \in [3,n-1] $. The `residual' polygon 
$ Z_{m},\ldots,Z_{n} $ can be rotated about the vertex 
$ Z_{m} = Z_{1} $ such that the `arriving' edgeline 
$ Z_{m-1}\vee Z_{m} $ and the `leaving' edgeline are
different and also that the triangle $ Z_{m-1},Z_{m},Z_{m+1} $
is negatively oriented. (The new vertices will not be
denoted anew.) This process doesn't change the sidelengths
nor the areas of the partial polygons 
$ Z_{1},\ldots,Z_{m-1} $ and $ Z_{m},\ldots,Z_{n} $ and of the 
whole polygon. Now, the vertex $ Z_{m} $ can be replaced by 
its mirror point $ Z_{m}' $ w.r.t. the line 
$ Z_{m-1}\vee Z_{m+1} $. This produces a positively oriented
kite quadrangle $ Z_{m-1},Z_{m}',Z_{m+1},Z_{m} $ as can be
seen from the standard position:

\vspace{-4ex}
\begin{center}
\begin{minipage}{8.3cm}
\bild{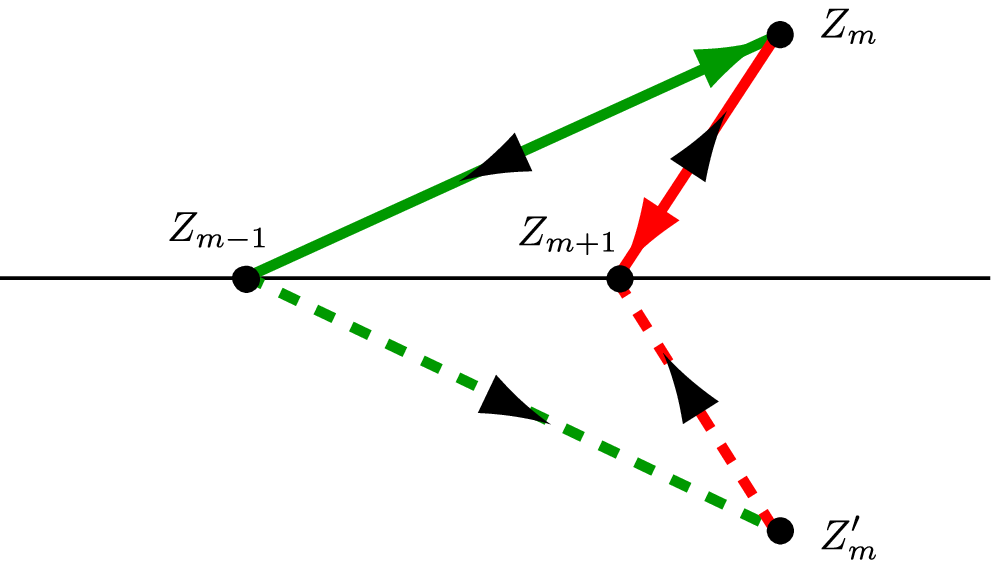}
\end{minipage}
\end{center}

\vspace{-4ex}
The kite quadrangle has positive area
$$
F' = 
\int_{[Z_{m-1},Z_{m}']} \omega +
\int_{[Z_{m}',Z_{m+1}]} \omega +
\int_{[Z_{m+1},Z_{m}]} \omega +
\int_{[Z_{m},Z_{m-1}]} \omega,
$$
and this implies that the new $ n $-gon is bigger
than the old one since
$$
\begin{aligned}
\int_{[Z_{m-1},Z_{m}']} \omega +
\int_{[Z_{m}',Z_{m+1}]} \omega
&=
F'-
\int_{[Z_{m+1},Z_{m}]} \omega -
\int_{[Z_{m},Z_{m-1}]} \omega \\
&=
F'+
\int_{[Z_{m},Z_{m+1}]} \omega +
\int_{[Z_{m-1},Z_{m}]} \omega \\
&>
\int_{[Z_{m-1},Z_{m}]} \omega +
\int_{[Z_{m},Z_{m+1}]} \omega.
\end{aligned}
$$
So, we can replace the polygon 
$ Z_{1},\ldots,Z_{m-1},Z_{m},Z_{m+1},\ldots,Z_{n} $ 
(in the new form modified by the above rotation) by the
bigger polygon 
$ Z_{1},\ldots,Z_{m-1},Z_{m}',Z_{m+1},\ldots,Z_{n} $ with
same sidelengths.
\QED

\medskip
As a result, a maximal $ n $-gon always has pairwise
distinct vertices.

\medskip
The next arguments concern the low cases $ n = 3,4 $.

\medskip
\BL
\label{lem_polymax_triangl}
If an oriented-convex triangle $ Z_{1}Z_{2}Z_{3} $ with
fixed sidelengths $ L_{1} $, $ L_{2} $ has maximal area then its 
parallelogram completion is cocyclic.
\EL

Of course, the third sidelength $ L_{3} $ is \TI{not} fixed here. 
The \TI{parallelogram completion} is the $ 4 $-gon 
$ Z_{1}Z_{2}Z_{3}Z_{4} $ such that $ Z_{4} $ is the point on
the opposite side of $ Z_{1} \vee Z_{3} $ to $ Z_{2} $ with 
$ d(Z_{3},Z_{4}) = L_{1} $ and $ d(Z_{1},Z_{4}) = L_{2} $.

\medskip
\TI{Proof of \ref{lem_polymax_triangl}.}
With fixed sidelengths $ L_{1} $, $ L_{2} $ and variable
enclosed angle $ \gamma $, the area of the triangle 
$ Z_{1}Z_{2}Z_{3} $ is given by the function 
$ f: \left]0,\pi\right[ \to \Rplus $:
$$
f(\gamma) := 
2\,\Arctan
\frac{\Sinh L_{1}\Sinh L_{2}\Sin \gamma}%
{(1+\Cosh L_{1})(1+\Cosh L_{2})-\Sinh L_{1}\Sinh
L_{2}\Cos\gamma};
$$
see Eqn. \eqref{eq_polyex_bil4}. The function $ f $ has the 
continuous extension $ f(0) := f(\pi) := 0 $. The derivative
can easily be calculated, and the condition $ f'(\gamma) = 0$
turns out to be equivalent to
\BE
\label{eq_polymax_triangl1}
\Cos \gamma = 
\frac{\Sinh L_{1}\Sinh L_{2}}
{(1+\Cosh L_{1})(1+\Cosh L_{2})}.
\EE
So, there is exactly one maximal point of $ f $, namely
contained in the open interval $ \left]0,\pi\right[ $.
Combined with the cosine law, it is seen that condition 
\eqref{eq_polymax_triangl1} is equivalent to 
\BE
\label{eq_polymax_triangl2}
\Cosh d(Z_{1},Z_{3}) = \Cosh L_{1}+ \Cosh L_{2}-1.
\EE
By the definition of $ Z_{4} $ and switching to the notation
of Sect. \ref{tripquad} we have 
$ K_{AE} = K_{BC} $, $ K_{CE} = K_{AB} $ and, from 
\eqref{eq_polymax_triangl2}, $ K_{AC} = K_{AB}+K_{BC}-1 $.
Substituting these values into the $ 17 $-identity 
(Lemma \ref{lem_tripquad_R17}) yields:
$$
\begin{aligned}
&(K_{AB}+K_{BC})\cdot(K_{BE}+1-K_{AB}-K_{BC})\cdot \\
&(-K_{BC}^2+K_{BC}K_{BE}+2K_{AB}K_{BC}-K_{BC}-2K_{BE}-K_{AB}^2
+2-K_{AB}+K_{BE}K_{AB}) = 0.
\end{aligned}
$$
The vanishing of the second and third parenthesis each time leads 
to a conditional equation for $ K_{BE} $ with the unique
solution
$$
K_{BE} = K_{AB}+K_{BC}-1 
\qquad \text{resp.} \qquad
K_{BE} = 
\frac{(K_{AB}-K_{BC})^2}{K_{AB}+K_{BC}-2}+1.
$$
This means, for the corresponding $ S $-values
\BE
\label{eq_polymax_triangl4}
S_{BE} = \sqrt{S_{AB}^2+S_{BC}^2}
\qquad \text{resp.} \qquad
S_{BE} =
\pm\,\frac{S_{AB}^2-S_{BC}^2}{\sqrt{S_{AB}^2+S_{BC}^2}}.
\EE
Transcribing Eqn. \eqref{eq_polymax_triangl2} to the 
$ S $-values gives
\BE
\label{eq_polymax_triangl3}
S_{AC} = \sqrt{S_{AB}^2+S_{BC}^2}.
\EE
Still, by assumption: $ S_{AE} = S_{BC} $ and 
$ S_{CE} = S_{AB} $. If this, together with 
\eqref{eq_polymax_triangl3} and the alternatives 
\eqref{eq_polymax_triangl4} is inserted into the Ptolemy
conditions \eqref{eq_tripquad_pt4} it turns 
out that always one of these conditions is satisfied, namely
the first one in case of the first alternative 
\eqref{eq_polymax_triangl4}, the second one in case of the 
plus sign, and the third one in case of the minus sign of
the second alternative in \eqref{eq_polymax_triangl4}.
\QED

\medskip
\BR
\label{rem_polymax_trianglopp}
The assumption that $ Z_{4} $ should lie on the other side
of $ Z_{2} $ w.r.t. $ Z_{1} \vee Z_{3} $ has not been used
in this proof. So, under the same hypothesis, also the points
$ Z_{1},Z_{2},Z_{3},Z_{4}^{\Ast} $ are cocyclic where 
$ Z_{4}^{\Ast} $ is the point on the same side of 
$ Z_{2} $ w.r.t. $ Z_{1} \vee Z_{3} $, satisfying
$ d(Z_{3},Z_{4}^{\Ast}) = L_{1} $ and 
$ d(Z_{1},Z_{4}^{\Ast}) = L_{2} $. However, this will not be
needed in the sequel.
\ER

The existence of polygons with prescribed sidelengths 
(part 1, Theorem 5.3) can be
completed to a continuous variant in case of triangles:

\medskip
\BL
\label{lem_polyex_contex}
Given positive real numbers $ L_{1},L_{2},L_{3} $ satisfying 
the three triangle inequalities in a strict manner, there
exist triangles in $ \B $ of a given orientation
with these sidelengths in continuous dependency on 
$ L_{1},L_{2},L_{3} $.
\EL

\Bew
For the vertices the following ansatz can be made:
$$
z_{1} =
\BM
1 \\
0 \\
0
\EM,
\qquad
z_{3} =
\BM
\Cosh L_{3} \\
\Sinh L_{3} \\
0
\EM,
\qquad
z_{3} =
\BM
\xi \\
\eta \\
\zeta
\EM,
\quad
\xi > 0, \quad
\xi^{2}-\eta^{2}-\zeta^{2} = 1.
$$
Then $ d(Z_{1},Z_{3}) = L_{3} $, and the demands 
$ d(Z_{1},Z_{2}) = L_{1} $, $ d(Z_{2},Z_{3}) = L_{2} $
translate to conditional equations for $ \xi,\eta,\zeta $
with the following solutions:
\BE
\label{eq_polyex_contex1}
\begin{aligned}
\xi &= \Cosh L_{1} \\[1ex]
\eta &= \frac{\Cosh L_{1}\Cosh L_{3}-\Cosh L_{2}}{\Sinh L_{3}} 
\\[1ex]
\zeta &= 
\frac{\pm2}{\Sinh L_{3}}
\sqrt{
\Sinh\frac{S}{2}
\Sinh\frac{D_{1}}{2}
\Sinh\frac{D_{2}}{2}
\Sinh\frac{D_{3}}{2}
}
\end{aligned}
\qquad \text{where} \qquad
\begin{aligned}
S &:= \phantom{-} L_{1}+L_{2}+L_{3} \\[1ex]
D_{1} &:= -L_{1}+L_{2}+L_{3} \\[1ex]
D_{2} &:= \phantom{-} L_{1}-L_{2}+L_{3} \\[1ex]
D_{3} &:= \phantom{-} L_{1}+L_{2}-L_{3}.
\end{aligned}
\EE
Here, the quantities $ D_{1},D_{2},D_{3} $ are just the
differences from the triangle quantities as in Eqn.
\eqref{eq_polymax_huil2}, so 
$ D_{1} \geq 0 $, $ D_{2} \geq 0 $, $ D_{3} \geq 0 $.
In front of the $ \zeta $ there occurs the minus resp. plus 
sign according as the triangle $ Z_{1},Z_{2},Z_{3} $ is
oriented positively or negatively. In particular, for
positive $ D_{k} $, these formulas show the
continuous (even real holomorphic) dependency as asserted.
\QED

\medskip
\BL
\label{lem_polymax_vier}
~\\[0.3ex]
If a  $ 4 $-gon $ Z_{1},Z_{2},Z_{3},Z_{4} $ in $ \B $ 
is maximal then it is cocyclic or non-strict (hence collinear).

\EL

\Bew
Consider in the strict case the following possibilities:

\smallskip
\TI{Case I.a:
Neither the points $ Z_{1},Z_{2},Z_{3} $ nor the points
$ Z_{3},Z_{4},Z_{1} $ are collinear.}

Both triangles are then strict. Look at further $ 4 $-gons 
$ Z_{1}',Z_{2}',Z_{3}',Z_{4}' $ with same sidelengths under
analogous assumptions.
These triangles are determined by the distance
$ L := d(Z_{1}',Z_{3}') $ if their orientations are kept
fixed. In order to apply differential calculus on the function 
$ L \mapsto F(Z_{1}',Z_{2}',Z_{3}')+ F(Z_{3}',Z_{4}',Z_{1}') $
it must be clarified that, for any $ L $ near 
$ L_{0} := d(Z_{1},Z_{3}) > 0 $, there are such further 
$ 4 $-gons. But this follows from Lemma 
\ref{lem_polyex_contex}, observing that the assumptions are
characterized by strict inequalities between continuous
functions. As a parameter for these neighbouring quadrangles
one may use the distance $ L $. (In the notation for these
quadrangles, the prime will be skipped from now on.) 

In addition, one may assume that both triangles 
$ Z_{1},Z_{2},Z_{3} $ and $ Z_{3},Z_{4},Z_{1} $ are positively
oriented since otherwise a genuine enlargement of the whole
area would be possible by reflections on the line 
$ Z_{1}\vee Z_{3} $. Then the points $ Z_{2} $ and $ Z_{4} $
lie on different sides of this line, and by Eqn.
\eqref{eq_polymax_bil_115} the whole area is given by 
$$
\begin{aligned}
& F(Z_{1},Z_{2},Z_{3},Z_{4}) =\\[1ex]
&=2\,\Arccos
\frac{\Cosh L_{1}+\Cosh L_{2}+\Cosh L+1}%
{\ds
4\,\Cosh\frac{L_{1}}{2}\Cosh\frac{L_{2}}{2}\Cosh\frac{L}{2}}
+
2\,\Arccos
\frac{\Cosh L_{3}+\Cosh L_{4}+\Cosh L+1}%
{\ds
4\,\Cosh\frac{L_{3}}{2}\Cosh\frac{L_{4}}{2}\Cosh\frac{L}{2}}.
\end{aligned}
$$
With the abbreviations
\BE
\label{eq_polymax_vier1}
a_{k} := \Cosh \frac{L_{k}}{2}, \quad k = 1,\ldots,4 \qquad
x := \Cosh \frac{L}{2}
\EE
and by means of the doubling formulas for the $ \cosh
$-function this amounts to the discussion of the function 
$$
f(x) := 
\Arccos \frac{a_{1}^2+a_{2}^2+x^2-1}{2a_{1}a_{2}x}
+
\Arccos \frac{a_{3}^2+a_{4}^2+x^2-1}{2a_{3}a_{4}x}
$$
or, by setting
$$
b_{12} := \frac{a_{1}^2+a_{2}^2-1}{2a_{1}a_{2}}, \qquad
c_{12} := \frac{1}{2a_{1}a_{2}}, \qquad
b_{34} := \frac{a_{3}^2+a_{4}^2-1}{2a_{3}a_{4}}, \qquad
c_{34} := \frac{1}{2a_{3}a_{4}},
$$
to
$$
f(x) = 
\Arccos\left(\frac{b_{12}}{x}+c_{12}x\right)+
\Arccos\left(\frac{b_{34}}{x}+c_{34}x\right).
$$
The derivative is calculated to be
$$
f' = A_{12}+A_{34},
$$
where
$$
A_{ij} := \frac{N_{ij}}{W_{ij}}, \quad
N_{ij} := \frac{b_{ij}}{x^2}-c_{ij}, \quad
W_{ij} := \sqrt{1-\left(\frac{b_{ij}}{x}+c_{ij}x\right)^2},
\qquad
ij \in \{12,34\}.
$$
The requirement $ f'(x) = 0 $ implies 
$ A_{12}^2-A_{12}^2 = 0 $. With regard to 
$ W_{ij} \neq 0 $ and $ x \neq 0 $ this can be written as an
algebraic equation $ P(x) = 0 $ where the polynomial 
$ P(x) $ is biquadratic (only containing $ x^4 $ and 
$ x^2 $). The search for its zeros is very much simplified
when one replaces the $ \cosh $-quantities 
\eqref{eq_polymax_vier1} by the corresponding 
$ \sinh $-quantities:
\BE
\label{eq_polymax_vier2}
\begin{gathered}
S_{AB} := \Sinh\frac{L_{1}}{2},\quad
S_{BC} := \Sinh\frac{L_{2}}{2},\quad
S_{CE} := \Sinh\frac{L_{3}}{2},\quad
S_{AE} := \Sinh\frac{L_{4}}{2}
\\
S_{AC} := \Sinh\frac{L}{2}.
\end{gathered}
\EE
The equation $ P(x) = 0 $ for $ x = \Cosh \frac{L}{2} $ is
then equivalent to an algebraic equation $ Q(y) = 0 $ for 
$ y = \Sinh \frac{L}{2} $. The decisive point is the
polynomial factorization $ Q(y) = Q_{+}(y) \cdot Q_{-}(y) $ 
where
$$
\begin{aligned}
Q_{+}(y) &:=
(S_{AB}S_{BC}+S_{CE}S_{AE})y^2-
(S_{AB}S_{AE}+S_{BC}S_{CE})(S_{AB}S_{CE}+S_{BC}S_{AE}) \\[1ex]
Q_{-}(y) &:=
(S_{AB}S_{BC}-S_{CE}S_{AE})y^2+
(S_{AB}S_{AE}-S_{BC}S_{CE})(S_{AB}S_{CE}-S_{BC}S_{AE}).
\end{aligned}
$$
If we had $ Q_{+}(S_{AC}) = 0 $ then we could
deduce the cocyclicity of the $ 4 $-gon $ ABCE $
from Theorem \ref{thm_tripquad_umkperron}.
If $ Q_{+}(S_{AC}) \neq 0 $ and thus $ Q_{-}(S_{AC}) = 0 $
nothing can be derived in the first instance, nevertheless in
connection with the next case:

\smallskip
\TI{Case I.b:
Neither the points $ Z_{2},Z_{3},Z_{4} $ nor the points
$ Z_{4},Z_{1},Z_{2} $ are collinear.} 

Here, completely analogous to case I.a, the 
maximality of the area implies 
$ R_{+}(S_{BE}) = 0 $ \TI{or} $ R_{-}(S_{BE}) = 0 $ for the two 
polynomials 
$$
\begin{aligned}
R_{+}(z) &:=
(S_{BC}S_{CE}+S_{AE}S_{AB})z^2-
(S_{BC}S_{AB}+S_{CE}S_{AE})(S_{BC}S_{AE}+S_{CE}S_{AB}) \\[1ex]
R_{-}(z) &:=
(S_{BC}S_{CE}-S_{AE}S_{AB})z^2+
(S_{BC}S_{AB}-S_{CE}S_{AE})(S_{BC}S_{AE}-S_{CE}S_{AB}).
\end{aligned}
$$
This case arises from case I.a by cyclically proceeding in
the list of vertices.

\smallskip
\TI{Case I:
Each three consecutive vertices of the $ 4 $-gon are not
collinear.}

Then both assumptions of the cases I.a and I.b are
satisfied, so $ Q_{+}(S_{AC})Q_{-}(S_{AC}) = 0 $ \TI{and} 
$ R_{+}(S_{BE})R_{-}(S_{BE}) = 0 $. If 
$ Q_{+}(S_{AC}) = 0 $ \TI{or} $ R_{+}(S_{BE}) = 0 $ the
cocyclicity is ensured by Theorem 
\ref{thm_tripquad_umkperron}. It remains the
discussion if $ Q_{-}(S_{AC}) = 0 $ \TI{and} 
$ R_{-}(S_{BE}) = 0 $. 

If both leading coefficients of the polynomials $ Q_{-} $, 
$ R_{-} $ don't vanish then the equations 
$ Q_{-}(S_{AC}) = R_{-}(S_{BE}) = 0 $ can be solved for 
$ S_{AC}^2 $ resp. $ S_{BE}^2 $, and the solutions imply
by multiplication
$$
S_{AC}^2 \cdot S_{BE}^2 =
(S_{AB}S_{CE}-S_{BC}S_{AE})^2.
$$
Thus one of the Ptolemy equations is satisfied, namely the
second or the third one in \eqref{eq_tripquad_pt4}.

If the leading coefficient of  $ Q_{-} $ is zero then so is 
the leading coefficient of $ R_{-} $, by the specific design
of $ R_{-} $. This yields $ S_{AB} = S_{AB} $ and 
$ S_{AE} = S_{BC} $. Then Lemma \ref{lem_polymax_triangl}
implies the cocyclicity of $ ABCE $ (by applying this lemma 
to the triangle $ ABC $).

If the leading coefficient of  $ R_{-} $ vanishes, so does 
the leading coefficient of  $ Q_{-} $ with the same result
as before.

\smallskip
\TI{Case II: 
The points $ Z_{1},Z_{2},Z_{3} $ are collinear and
$ Z_{4} $ is not on the line of $ Z_{1},Z_{2},Z_{3} $.}

Set $ D_{13} := d(Z_{1},Z_{3}) > 0 $. Without loss, we may
assume $ F(Z_{3},Z_{4},Z_{1}) > 0 $. 

Three collinear and pairwise distinct points always form a
non-strict $ 3 $-gon. For a non-strict $ 3 $-gon, exactly
one of the triangle inequalities is converted to an equality.
So there exist three subcases
\begin{align}
\label{eq_polymax_V1}
\tag{II.a}
&D_{13} = L_{1}+L_{2}, \qquad 
L_{1} < L_{2}+D_{13}, \quad
L_{2} < L_{1}+D_{13} \\[1ex]
\label{eq_polymax_V2}
\tag{II.b}
&D_{13} < L_{1}+L_{2}, \qquad 
L_{1} = L_{2}+D_{13}, \quad
L_{2} < L_{1}+D_{13} \\[1ex]
\label{eq_polymax_V3}
\tag{II.c}
&D_{13} < L_{1}+L_{2}, \qquad 
L_{1} < L_{2}+D_{13}, \quad
L_{2} = L_{1}+D_{13}.
\end{align}
In all these cases, a suitable variation of $ D_{13} $
produces a $ 4 $-gon with larger area and like sidelengths 
$ L_{1},L_{2},L_{3}, L_{4} $. This is possible since the
area of $ Z_{3}Z_{4}Z_{1} $ depends differentiably on 
$ D_{13} $ while for the area of $ Z_{1}Z_{2}Z_{3} $ one of 
the relations \eqref{eq_polymax_huil4}, 
\eqref{eq_polymax_huil3} is appropriate. 
For example, in
case \eqref{eq_polymax_V1} one has to diminish $ D_{13} $ in
such a way that in all relations of this row the
genuine smaller sign remains. Then, for the limit 
$ D_{13} \upto L_{1}+L_{2} $, Eqn. \eqref{eq_polymax_huil4}
is effective such that the area development of 
$ Z_{1},Z_{2},Z_{3} $ exceeds that of $ Z_{3},Z_{4},Z_{1} $ 
if $ D_{13} $ is sufficiently close to $ L_{1}+L_{2} $.
Similar arguments apply in the cases 
\eqref{eq_polymax_V2} and \eqref{eq_polymax_V3}, where here 
$ D_{13} $ must be enlarged suitably, and Eqn. 
\eqref{eq_polymax_huil3} must be used instead.

So in case II the given quadrangle cannot be maximal.

\smallskip
\TI{Case III: 
The vertices $ Z_{1},Z_{2},Z_{3} $ are collinear as well as 
the vertices $ Z_{3},Z_{4},Z_{1} $.}

Then all vertices are collinear (by $ Z_{1} \neq Z_{3} $)
and the whole area is zero. On the other hand there exists a
cocyclic and oriented-convex $ 4 $-gon with same sidelengths
by part 1, Theorem 5.3. This $ 4 $-gon has positive
area (Lemma \ref{lem_polymax_nn}) such that the original area
was not maximal.
\QED

\medskip
\TI{Remark.}
The special case of Lemma \ref{lem_polymax_vier} for three
\TI{equal} sidelengths has been discussed in Leichtwei"s
[2005], Lemma 5.9 with a different method. Generally, in
this paper, Leichtwei"s solved the maximum area problem for
curves in $ \B $ of fixed constant width.

\medskip
\BC
\label{cor_polymax_neck}
If a $ n $-gon $ Z_{1} \ldots Z_{n} $ in $ \B $ has 
maximal area compared to all $ n $-gons with same
sidelengths then it is either cocyclic or else non-strict 
(hence collinear).
\EC

\Bew
Consider the following cases I and II:

\smallskip
\TI{Case I: There exists an index $ k \in [3,n] $ such that 
the polygon $ Z_{1} \ldots Z_{k} $ is not strict.}
Then choose $ k $ \TI{maximal with this property.}

If $ k = n $ we are done.

Now assume $ k \leq n-1 $. Then $ Z_{1},\ldots,Z_{k} $ are
collinear, due to part 1, Theorem 5.3 (ii). By the maximal 
choice of $ k $, the polygon $ Z_{1} \ldots Z_{k+1} $ is
strict. The point $ Z_{k+1} $ is then outside the line of 
$ Z_{1},\ldots,Z_{k} $ because otherwise the area of the
polygon $ Z_{1},\ldots,Z_{k+1} $ could be increased (from 
$ 0 $ to a positive value), according to 
part 1, Theorem 5.3 (i). Also the polygon 
$ Z_{1},\ldots,Z_{n} $ could be enlarged. By the same
reason, the $ 4 $-gon $ Z_{k-2}Z_{k-1}Z_{k}Z_{k+1} $ cannot 
be enlarged (it cannot be collinear). But a hyperbolic circle
cannot contain three pairwise distinct collinear points.
Thus this situation is not existent.

\smallskip
\TI{Case II: For any index $ k \in [3,n] $ the polygon
$ Z_{1} \ldots Z_{k} $ is strict.}

The triangle $ Z_{1},Z_{2},Z_{3} $ cannot be enlarged
because otherwise the whole $ n $-gon could be made larger.
This triangle is not collinear since otherwise a
non-collinear triangle with positive area and same
sidelengths would exist (part 1, Theorem 5.3 (i)).
Thus, the $ 3 $-gon $ Z_{1}Z_{2}Z_{3} $ is cocyclic, as any
non-collinear $ 3 $-gon.

Now, an induction can be started:

The $ 4 $-gon $ Z_{1}Z_{2}Z_{3}Z_{4} $ is strict and cannot 
be enlarged, hence it is cocyclic (Lemma \ref{lem_polymax_vier}).

The $ 4 $-gon $ Z_{2}Z_{3}Z_{4}Z_{5} $ is strict. Otherwise,
these points were collinear what is not possible with regard
to the points $ Z_{2},Z_{3},Z_{4} $; these points are
cocyclic by the above argument. Moreover, the two
circum-circles of $ Z_{1},Z_{2},Z_{3},Z_{4} $ 
and of $ Z_{2},Z_{3},Z_{4},Z_{5} $ are equal since they have
the points $ Z_{2},Z_{3},Z_{4} $ in common.

Obviously, one can proceed in the same way and realize 
successively that the $ 4 $-gons 
$$ 
Z_{3}Z_{4}Z_{5}Z_{6}, \qquad Z_{4}Z_{5}Z_{6}Z_{7}, 
\qquad \ldots 
$$
are all cocyclic with a fixed circum-circle. Thus, the whole
$ n $-gon is cocyclic.
\QED

\medskip
In order to establish the oriented convexity of the maximal 
polygons we use a reduction lemma which is well known in
Euclidean convexity (see e.g. \TI{Moret/Shapiro} [1990] or 
\TI{Pinelis} [2006]). In the present context it reads as:

\medskip
\BL[reduction lemma]
\label{lem_polyrig_ind}
For $ n \geq 3 $ let points $ Z_{1},\ldots,Z_{n+1} \in \B $ 
be given. If, for any $ k \in \{1,\ldots,n+1\} $, the points
$ Z_{1},\ldots,\Hat{Z_{k}},\ldots, Z_{n+1} $ form an
oriented-convex $ n $-gon then the points 
$ Z_{1},\ldots,Z_{n+1} $ form an oriented-convex $ (n+1) $-gon.
\EL

The roof over an element in a list means omission of the
element.

\medskip
\TI{Proof of \ref{lem_polyrig_ind}.}
The argument is not very different from the Euclidean
situation since the convexity notions are rather near in
both geometries. 

Let the given $ n+1 $ points be denoted somewhat differently,
namely as $ Z_{1},\ldots,Z_{n},Z^{\Ast} $. By hypothesis, the polygon 
$ Z_{1},\ldots,Z_{n} $ is oriented-convex, and for \TI{this}
polygon the cyclic index convention from part 1 will be
maintained, so $ Z_{n+1} := Z_{1} $, etc.

In the polygon $ Z_{1} \ldots Z_{n}Z^{\Ast} $ there is one 
additional point compared to the polygon $ Z_{1} \ldots Z_{n} $ and 
there are two edgelines more, namely 
$ \pfeil{Z_{n} \vee Z^{\Ast}} $ and 
$ \pfeil{Z^{\Ast} \vee Z_{1}} $. (In return, the edgeline 
$ \pfeil{Z_{n} \vee Z_{1}} $ is omitted.) In order to gain
the defining conditions of oriented convexity 
(part 1, Sect. 3) for $ Z_{1},\ldots,Z_{n},Z^{\Ast} $ we must show:
\begin{DES1}{\textbullet}
\item[\textbullet]
the relations of the `new' point $ Z^{\Ast} $ with the `old'
edgelines, i.e.
\BE
[Z_{k},Z_{k+1},Z^{\Ast}] > 0, \quad k = 1,\ldots,n-1;
\label{eq_polyrig_ind1}
\tag{a}
\EE
\vspace{-4ex}
\item[\textbullet]
the relations of the two `new' edgelines with the `old'
points, i.e.
\begin{alignat}{2}
[Z_{n},Z^{\Ast}, Z_{k}] &> 0, &\quad & k = 1,\ldots,n-1
\label{eq_polyrig_ind2}
\tag{b}
\\[0.5ex]
[Z^{\Ast}, Z_{1},  Z_{k}] &> 0, &\quad & k = 2,\ldots,n.
\label{eq_polyrig_ind3}
\tag{c}
\end{alignat}
\end{DES1}
\TI{For \eqref{eq_polyrig_ind1}:}
Among the old vertices one can omit \TI{one}, e.g. one
with an index $ j \in \{1,\ldots,k-1\} $ (only possible for 
$ k \geq 2 $) or one with an index $ j \in \{k+2,\ldots,n\} $ 
(only possible for $ k \leq n-2 $). In the first case, the
relation $ [Z_{k},Z_{k+1},Z^{\Ast}] > 0 $ is deduced from
the oriented convexity of 
$$
Z_{1} \ldots \Hat{Z_{j}} \ldots Z_{k}Z_{k+1} \ldots 
Z_{n}Z^{\Ast},
$$
and in the second case from the oriented convexity of 
$$
Z_{1} \ldots Z_{k}Z_{k+1}
\ldots \Hat{Z_{j}} \ldots Z_{n} Z^{\Ast}.
$$
\TI{For \eqref{eq_polyrig_ind2}:}
Here, one succeeds by omitting an index 
$ j \in \{1,\ldots,n\} \setminus \{k,n\} $, meaning that one
considers the oriented-convex polygon 
$ Z_{1} \ldots \Hat{Z_{j}} \ldots Z_{n}Z^{\Ast} $ and
deducing from this $ [Z_{n},Z^{\Ast}, Z_{k}] > 0 $.

\TI{For \eqref{eq_polyrig_ind3}:}
Again, by omitting an index 
$ j \in \{1,\ldots,n\} \setminus \{1,k\} $, one deduces 
from the oriented convexity of 
$ Z_{1} \ldots \Hat{Z_{j}} \ldots Z_{n}Z^{\Ast} $ that
$ [Z^{\Ast}, Z_{1},  Z_{k}] > 0 $.
\QED

\medskip
\BL
\label{lem_polymax_okonv}
If a cocyclic $ n $-gon $ Z_{1} \ldots Z_{n} $ in $ \B $ has 
maximal area compared to all $ n $-gons with same
sidelengths then it is oriented-convex.
\EL

\Bew
From the maximality follows at any rate that the vertices
are pairwise distinct and from the cocyclicity that the 
$ n $-gon is not collinear. 
Now the proof proceeds by induction on $ n $, using the
foregoing Lemma \ref{lem_polyrig_ind}.

\TI{Initial step $ n = 3 $:}
For three non-collinear points $ Z_{1},Z_{2},Z_{3} $ there
are only the two possibilities $ [Z_{1},Z_{2},Z_{3}] < 0 $
and $ [Z_{1},Z_{2},Z_{3}] > 0 $. Due to the maximality, only
the second possibility is left over. 

\TI{Induction step from $ n $ to $ n+1 $ for $ n = 3 $:}
Let the $ (n+1) $-gon $ Z_{1} \ldots Z_{n+1} $ be cocyclic
and of maximal area. Consider, for any 
$ k \in \{1,\ldots,n+1\} $, the points 
$ Z_{1},\ldots,\Hat{Z_{k}},\ldots,Z_{n+1} $. They form a
cocyclic $ n $-gon (by $ Z_{k-1} \neq Z_{k+1} $). It has 
maximal area because otherwise the polygon 
$ Z_{1} \ldots Z_{n+1} $ could be increased. By the
induction hypothesis, the  $ n $-gon 
$ Z_{1} \ldots \Hat{Z_{k}} \ldots Z_{n+1} $ is
oriented-convex. This being true for any $ k $, the
preceding reduction Lemma \ref{lem_polyrig_ind} shows that
also the polygon $ Z_{1},\ldots,Z_{n+1} $ must be
oriented-convex.
\QED

\medskip
The main theorem now arises by combining the above results
with those of part 1:

\medskip
\BT
\label{thm_polymax_main}
For any $ n $-gon in the hyperbolic plane with sidelengths 
$ L_{1},\ldots,L_{n} $ there exist a $ n $-gon of maximal
area with same sidelengths. This maximal copy is either
cocyclic and oriented-convex or else collinear and
monotonically arranged.

In both cases the maximal copy is uniquely determined up to 
hyperbolic motion. In the first case, its area is positive, 
in the second case the area vanishes. 

Which case occurs is solely determined by the behaviour of
the sidelengths as real numbers: If, for all 
$ k = 1,\ldots,n $, there holds
$$
L_{k} < L_{1}+\cdots+\Hat{L_{k}}+\cdots+L_{n}
$$
then the first case is present, otherwise the second. 
\ET

\Bew

The \TI{existence} follows from Lemma
\ref{lem_polymax_genmax}, the alternative from Corollary 
 \ref{cor_polymax_neck}, the oriented convexity from Lemma 
\ref{lem_polymax_okonv}, the sign of the maximal area from
Lemma \ref{lem_polymax_nn}, and the
uniqueness from part 1, Theorem 5.3.
\QED

\newpage
\bigskip
\TB{\Large References}

\smallskip
Bilinski, S. [1969]:
Zur Begr"undung der elementaren Inhaltslehre
in der hyperbolischen Ebene:
Math. Ann. 180, 256-268

Blaschke, W. [1956]:
Kreis und Kugel:
2. durchgesehene und verbesserte Auflage, 
Berlin: Walter de Gruyter 1956

Blumenthal, L.M. [1970]:
Theory and applications of distance geometry:
Second Edition:
Chelsea Publ. Comp. Bronx, New York, i-xi and 1-347

Cartan, H. [1967]:
Formes diff\'erentielles:
Hermann \& Cie. Paris, 1-185

Knebelman, M.S. [1941]:
Two isoperimetric problems:
The American Mathematical Monthly 48,
623-627

Kryzhanovskij, D.A. [1959]:
Izoperimetry (Russian):
Moskau: Gosdarstvennoe izdatelstvo
fisiko-matematiceskoj literatury 

Leichwei"s, K. [2005]:
Curves of constant width in the non-Euclidean geometry:
Abh. Math. Semin. Univ. Hamb. 75, 257-284

Menger, K. [1928]:
Untersuchungen "uber allgemeine Metrik:
Math. Ann. 100, 75-163

Moret, B.M.E. and Shapiro, H.D. [1990]:
Algorithms from P to NP. Vol. I: design and efficiency:
Addison-Wesley Amsterdam, I- XV, 1-576 

Perron, O. [1964]:
Seiten und Diagonalen eines Kreisvierecks
in der hyperbolischen Geometrie:
Math. Z. 84, 88-92

Pinelis, I. [2006]:
Convexity of sub-polygons of convex polygons:
arXiv:math/0609698v1, 1-24

van der Waerden, B.L. [2003]
Algebra, Volume I
Springer Berlin, i-xiv, 1-265
(English translation of the 7th ed.)

Walter, R. [2010]:
Polygons in hyperbolic geometry 1:
Rigidity and inversion of the $ n $-inequality:
arXiv:1008.3404v1

Yaglom, I.M. und Boltjanski, W.G. [1951]:
Konvexe Figuren:
Berlin: VEB Deutscher Verlag der Wissenschaften
(German translation of the Russian original:
Bibl. Math. Zirk. 4. Moskau-Leningrad: Staatsverlag technisch- 
theoretische Literatur 1951) 

\medskip
Rolf Walter\\
Fakult"at f"ur Mathematik\\
Technische Universit"at Dortmund\\
Arbeitsgebiet Differentialgeometrie\\
Vogelpothsweg 87\\
D-44227 Dortmund

E-Mail: rolf.walter$ @ $tu-dortmund.de\\
URL: http://www.mathematik.uni-dortmund.de/$\sim$walter/

\end{document}